\documentclass[11pt]{article} \textheight = 24truecm \textwidth
= 16truecm \hoffset = -2truecm \voffset = -2truecm
\usepackage[english, russian]{babel}
\usepackage{amssymb,amsthm}
\usepackage{amsmath}

\begin{document}

\large \sloppy

\begin{center}
\textsc {On uniformly metrizability of the functor of idempotent
probability measures}

\textsl{A. A. Zaitov, I. I. Tojiev}
\end{center}

\begin{abstract}
In the present paper we show that the functor of idempotent
probability  measures satisfies all of conditions with an
additional claim of uniform metrizability of functors.
\end{abstract}



The present paper is a continuation of [1]. We begin it with some
definitions from [2].

\textbf{Definition 1}. A functor $\mathcal{F}$ acting in the
category $\textit{Comp}$ of Hausdorff compact spaces and their
continuous mappings is called to be \textit{seminormal} if it
satisfies the following conditions:

\textit{1) $\mathcal{F}$ preserves empty set and singleton}, i. e.
$\mathcal{F}(\emptyset)=\emptyset$ and $\mathcal{F}(1)=1$ take
place, where $1$ is a singleton.

\textit{2) $\mathcal{F}$ preserves intersections}, i. e. for  a
given compacta $X$ and for every family $\mathcal{B}$ of closed
subsets of $X$ the equality $\mathcal{F}\left(\bigcap\limits_{F\in
\mathcal{B}}F\right)=\left(\bigcap\limits_{F\in
\mathcal{B}}\mathcal{F}(F)\right)$ holds;

\textit{3) $\mathcal{F}$ is monomorphic}, i. e. for any embedding
$i:A\rightarrow X$ the map
$\mathcal{F}(i):\mathcal{F}(A)\rightarrow \mathcal{F}(X)$ is also
embedding;

\textit{4) $\mathcal{F}$ is continuous}, i. e. for any spectrum
$S=\{X_\alpha,\ \pi_\alpha^\beta;\ \mathrm{A}\}$ we have
$\mathcal{F}(\lim S)=\lim(\mathcal{F}(S))$.

If a functor $\mathcal{F}$ is seminormal then there exists unique
natural transformation
$\eta^{\mathcal{F}}=\eta:Id\rightarrow\mathcal{F}$ of identity
functor $Id$ into functor $\mathcal{F}$. Moreover this
transformation is monomorphism, i. e. for each Hausdorff compact
space $X$ the map $\eta^{\mathcal{F}}:X\rightarrow\mathcal{F}(X)$
is embedding.

\textbf{Definition 2.} A seminormal functor $\mathcal{F}$, acting
in the category $\textit{MComp}$ of metrizable compact spaces is
called to be \textit{metrizable} if for any metrizable compact $X$
and for each metric $d=d_X$ on $X$ it is possible to put a
conformity the metric $d_{\mathcal{F}(X)}$ on compact
$\mathcal{F}(X)$ such that the following conditions hold:

Ð1) if $i:(X_1, d^1)\rightarrow (X_2, d^2)$ is isometrical
embedding then $\mathcal{F}(i): (\mathcal{F}(X_1),
d_{\mathcal{F}(X_1)}^1)\rightarrow (\mathcal{F}(X_2),
d_{\mathcal{F}(X_2)}^2)$ is also isometrical embedding;

Ð2) the embedding $\eta_X:(X, d)\rightarrow(\mathcal{F}(X),
d_{\mathcal{F}(X)})$  is isometric;

Ð3) $diam\mathcal{F}(X)=diamX$ .

\textbf{Definition 3.} A metrizable functor $\mathcal{F}$ is
called to be \textit{uniform metrizable}, if its some metrication
has the property

Ð4) for any continuous mapping $f:(X_1, d^1)\rightarrow(X_2,
d^2)$  the mapping $\mathcal{F}^+(f): (\mathcal{F}^+(X_1),
d_+^1)\rightarrow(\mathcal{F}^+(X_2), d_+^2)$ is uniform
continuous\footnote{For definition of $\mathcal{F}^+$ in case of the
functor of idempotent probability measures, see below.}.\\

Let $S$ be a set equipped with two algebraic operation: addition
$\oplus$ and multiplication $\odot$. $S$ is called [3] a semiring
if the following conditions hold:

($i$) the addition $\oplus$ and the multiplication $\odot$ are
associative;

($ii$) the addition $\oplus$ is commutative;

($iii$) the multiplication $\odot$ is distributive with respect to
the addition $\oplus$.

A semiring $S$ is commutative if the multiplication $\odot$ is
commutative. A unity of semiring $S$ is an element $\textbf{1}\in
S$ such that $\textbf{1}\odot x=x\odot\textbf{1}=x$ for all $x\in
S$. A zero of a semiring $S$ is an element $\textbf{0}\in S$ such
that $\textbf {0}\neq \textbf{1}$ and $\textbf{0}\oplus x=x$,
$\textbf{0}\odot x=x\odot \textbf{0}=\textbf{0}$ for all $x\in S$.
A semiring $S$ is idempotent if $x\oplus x=x$ for all $x\in S$. A
semiring $S$ with zero $\textbf{0}$ and unity $\textbf{1}$ is
called a semifield if each nonzero element $x\in S$ is invertible.

Let $\mathbb{R}$ be the field of real numbers and $\mathbb{R_{+}}$
the semifield of nonnegative real numbers (with respect to the
usual operations). The change of variables $x\mapsto u=h\ln x$,
$h>0$, defines a map $\Phi_{h}:\mathbb{R_{+}}\rightarrow
S=\mathbb{R}\cup\{-\infty\}$. Let the operations of addition
$\oplus$ and multiplication $\odot$ on $S$ be the images of the
usual operations of addition $+$ and multiplication $\cdot$ on
$\mathbb{R}$, respectively, by the map $\Phi_{h}$, i. e. let
$u\oplus_{h}v=h \ln(\exp (u/h)+\exp(v/h)),$ $u\odot v=u+v$. Then
we have $\textbf{0}=-\infty=\Phi_{h}(0)$,
$\textbf{1}=0=\Phi_{h}(1)$. It is easy to see that
$u\oplus_{h}v\rightarrow \max\{u, v\}$ as $h\rightarrow 0$. Hence,
$S$ forms semifield with respect to operations $u\oplus
v=\textrm{max}\{u,v\}$ and $u\odot v=u+v$. It denotes by
$\mathbb{R}_{\textrm{max}}$. It is idempotent. This passage from
$\mathbb{R}_+$ to $\mathbb{R}_{\textrm{max}}$ is called the Maslov
dequantization.

Let $X$ be a compact Hausdorff space, $C(X)$ the algebra of
continuous functions $\varphi : X \rightarrow \mathbb{R}$ with the
usual algebraic operations. On $C(X)$ the operations $\oplus$ and
$\odot$ define as follow:

$\varphi \oplus \psi=\textrm{max}\{\varphi, \psi \}$, where
$\varphi, \psi \in C(X)$,

$\varphi \odot \psi=\varphi + \psi$, where $\varphi$, $\psi \in
C(X)$,

$\lambda \odot \varphi=\varphi+\lambda_{X}$, where $\varphi\in
C(X)$, $\lambda\in \mathbb{R},$  and $\lambda_X$ is a constant
function.

Recall [4] that a functional $\mu : C(X)\rightarrow
\mathbb{R}(\subset \mathbb{R}_{\textrm{max}})$ is called to be an
idempotent probability measure on $X$, if:

1)  $\mu (\lambda_{X})=\lambda$ for each $\lambda \in \mathbb{R}$;

2)  $\mu (\lambda \odot \varphi)=\mu (\varphi)+\lambda$ for all
$\lambda \in \mathbb{R}$, $\varphi \in C(X)$;

3)  $\mu (\varphi \oplus \psi)=\mu(\varphi)\oplus \mu(\psi)$ for
every $\varphi$, $\psi\in C(X)$.

The number $\mu(\varphi)$ is named Maslov integral of $\varphi\in
C(X)$ with respect to $\mu$.

For a compact Hausdorff space $X$ a set of all idempotent
probability measures on $X$ denotes by $I(X)$.  Consider $I(X)$ as
a subspace of $\mathbb{R}^{C(X)}$. In the induced topology the
sets
\begin{center} $\langle \mu; \varphi_1, \varphi_2, ...,
\varphi_k; \varepsilon \rangle=\{\nu\in I(X):
|\mu(\varphi_i)-\nu(\varphi_i)|<\varepsilon, i=1, ..., k\}$,
\end{center}
form a base of neighborhoods of the idempotent measure $\mu\in
I(X)$, where $\varphi_i\in C(X)$, $i=1, ..., k$, and $\varepsilon
>0$. The topology generated by this base coincide with pointwise
topology on $I(X)$. The topological space $I(X)$ is compact [4].
Given a map $f:X\rightarrow Y$ of compact Hausdorff spaces the map
$I(f):I(X)\rightarrow I(Y)$ defines by the formula
$I(f)(\mu)(\varphi)=\mu(\varphi\circ f)$, $\mu\in I(X)$, where
$\varphi\in C(Y)$. Thus the construction $I$ is a covariant
functor, acting in the category of compact Hausdorff spaces and
their continuous mappings. As it is known [4] the functor is
normal in Schepin's sense, let us check if it is metrizable.

For any given idempotent measure $\mu\in I(X)$ we may define the
support of $\mu$:
\begin{center}supp $\mu=\bigcap\{A\subset X:
\overline{A}=A,$ $\mu\in I(A)\}$.
\end{center}

Let $\rho :X\times X\rightarrow \mathbb{R}$ be a metric, and
$\rho_I :I(X)\times I(X)\rightarrow \mathbb{R}$ be as in
[1]\footnote{The secondary author calls $\rho_I$ as 'Zaitov
metric'.}.

\textbf{Lemma 1.} \textsl{Let $X$ be a metric space with metric
$\rho$. Then $\delta_X:(X,\rho)\rightarrow (I(X), \rho_I)$ is an
isometry.}

\textsc{Proof}. For any pair $x_1, x_2\in X$ one has
$\delta_{x_1}, \delta_{x_2}\in I(X)$, and
$$
\rho_I(\delta_{x_1},\  \delta_{x_2})=\rho_{\omega}(\delta_{x_1},\
\delta_{x_2})= \rho_{\omega}(0\odot\delta_{x_1},\
0\odot\delta_{x_2})=
$$
$$
=\text{min}\left \{ diamX,\bigoplus\limits_{(x_1,\ x_2)\in S\xi}
|0-0|\odot\rho(x_1,\ x_2)\right\}=
$$
$$
=\text{min}\{diamX,\rho(x_1,\ x_2)\}= \rho(x_1,\ x_2).
$$

Lemma 1 is proved.

\textbf{Lemma 2.} \textsl{For any metric on the compactum $X$ the
following equality holds}
$$diam(X,\ \rho)=diam(I(X),\ \rho_I).$$

\textsc{Proof.} Identify each point $x\in X$ with Dirac measure
$\delta_x\in I(X)$, which gives embedding $X\subset_{\rightarrow}
I(X)$. Hence by Lemma 1 one has $diamX\leq diamI(X)$. Now we show
$diamI(X) \leq diamX$. Let $\mu_k\in I(X)$, $k=1,2$, be an
arbitrary pairs of idempotent measures. Consider sequences
$\{\mu_k^{(n)}\}_{n=1}^{\infty}\subset I_{\omega}(X)$, $k=1,\ 2$,
such that $\mu_k^{(n)}\rightarrow\mu_k$. Then according to
definition of $\rho_I$ (see formula (6) [1]) we have
$\rho_I(\mu_1,\mu_2)=\lim\limits_{n\rightarrow\infty}\rho_{\omega}(\mu_1^{(n)},\mu_2^{(n)})$.
The definition of $\rho_\omega$ for all $\mu_1^{(n)},\
\mu_2^{(n)}\in I_{\omega}(X)$ implies the following inequality
$$
\rho_{\omega}(\mu_1^{(n)},\mu_2^{(n)})=\text{min}\left\{diamX,
\bigoplus\limits_{(x_{1j},x_{2k})\in
S\xi}|\lambda_{1j}-\lambda_{2k}|\odot\rho(x_{1j},x_{2k})\right\}\leq
diamX.
$$
From here one has
$\rho_I(\mu_1,\mu_2)=\lim\limits_{n\rightarrow\infty}\rho_{\omega}(\mu_1^{(n)},\mu_2^{(n)})\leq
diamX$, and by forcing of arbitrariness of $\mu_1,\ \mu_2\in I(X)$
it follows $diamI(X)\leq diamX$. Lemma 2 is proved.

\textbf{Lemma 3.} \textsl{Let $(X_1,\rho^1)$, $(X_2,\rho^2)$ be
metrizable compacta such that $diam(X_1,\rho^1)=diam(X_2,\rho^2)$.
If $i: (X_1,\rho^1)\rightarrow (X_2,\rho^2)$ is an isometrical
embedding then $I(i):(I(X_1),\rho^1_{I,\
X_1})\rightarrow(I(X_2),\rho^2_{I,\ X_2})$ is also isometrical
embedding.}

\textsc{Proof.} Note that the condition
$diam(X_1,\rho^1)=diam(X_2,\rho^2)$ in Lemma 3 is essentially.
Really let $(X_1,\rho^1)$, $(X_2,\rho^2)$ be metric spaces and
what's more $diam(X_1,\rho^1)<diam(X_2,\rho^2)$, and let $\zeta
:X_1\rightarrow X_2$ be an isometrical embedding. Take arbitrary
points $x_1, x_2\in X_1$. Consider non-positive number
$\lambda_1,\lambda_2\in \mathbb{R}_{\text{max}}$ such that
$diam(X_2,\rho^2)<|\lambda_1-\lambda_2|$. For the idempotent
probability measures
$$
\mu_1=0\odot\delta_{x_1}\oplus\lambda_1\odot\delta_{x_2}
$$
and
$$
\mu_2=0\odot\delta_{x_1}\oplus\lambda_2\odot\delta_{x_2}
$$
it is clear that $supp\mu_1=supp\mu_2=\{x_1, x_2\}$. Hence by the
definition
$$
\rho_{\omega}^{X_1}(\mu_1,\mu_2)=\text{min}\{diam(X_1,\rho^1),
|\lambda_1-\lambda_2|\}=diam(X_1,\rho^1).
$$

Repeating this procedure for the idempotent probability measures
$I(i)(\mu_1)$ and $I(i)(\mu_2)$ we get
$$
\rho_{\omega}^{X_2}(I(i)(\mu_1), I(i)(\mu_2))=diam(X_2,\rho^2)
$$

Thus
$\rho_{\omega}^{X_1}(\mu_1,\mu_2)\neq\rho_{\omega}^{X_2}(I(i)(\mu_1),I(i)(\mu_2))$.

Let now we have $diam(X_1,\rho^1)=diam(X_2,\rho^2)$. By the
definition of $\rho_I$ it is enough to consider idempotent
probability measures
$\mu_{k}=\lambda_{k1}\odot\delta(x_{k1})\oplus...\oplus\lambda_{kn_k}\odot\delta(x_{kn_k})$,
$k=1\ 2$. Then by the definition we have
$$
I(i)(\mu_{k})(\varphi)=\mu_{k}(\varphi\circ i)=
(\lambda_{k1}\odot\delta(x_{k1})\oplus...\oplus\lambda_{kn_k}\odot\delta(x_{kn_k}))(\varphi\circ
i)=
$$
$$
=\lambda_{k1}\odot(\delta(x_{k1})(\varphi\circ
i))\oplus...\oplus\lambda_{kn_k}\odot(\delta(x_{kn_k})(\varphi\circ
i))=\lambda_{k1}\odot\varphi(i(x_{k1}))\oplus...\oplus\lambda_{kn_k}\odot\varphi(i(x_{kn_k}))=
$$
$$
=\lambda_{k1}\odot\delta(i(x_{k1}))(\varphi)\oplus...\oplus\lambda_{kn_k}\odot\delta(i(x_{kn_k}))(\varphi)
=(\lambda_{k1}\odot\delta(i(x_{k1}))\oplus...\oplus\lambda_{kn_k}\odot\delta(i(x_{kn_k})))(\varphi),
$$
i. e.
$I(i)(\mu_{k})=\lambda_{k1}\odot\delta(i(x_{k1}))\oplus...\oplus\lambda_{kn_k}\odot\delta(i(x_{kn_k}))$.
That is why $\rho^2_{I,\ X_2}(I(i)(\mu_1),\
I(i)(\mu_2))=\rho^1_{I,\ X_1}(\mu_1,\ \mu_2)$. Lemma 3 is proved.

Let now we show that the functor $I$ satisfies property Ð4) with
an  additional condition, more exactly with condition of equality
of diameters of consider compacta. For this we need the following
construction. Since functor $I$ is normal there exists unique
natural transformation $\eta^I=\eta:Id\rightarrow I$ of identity
functor $Id$ into functor $I$. Here the natural transformation
$\eta$ consists of monomorphisms $\delta_X$, $X\in\textit{Comp}$.
More detail the last means that for each compact $X$ the mapping
$\delta_X:X\rightarrow I(X)$, which defines as
$\delta_X(x)=\delta_x$, $x\in X$, is an embedding. Thus
$\eta=\{\delta_X:X\in\textit{Comp}\}$.

Let $X$ be a metrizable compact. Put $I^0(X)=X$,
$I^k(X)=I(I^{k-1}(X))$, $k=1,2,...$ and
$\eta_{n-1,n}=\eta_{I^{n-1}(X)}:I^{n-1}(X)\rightarrow I^n(X)$. For
$n<m$ denote
$$
\eta_{n,m}=\eta_{m-1,m}\circ...\circ\eta_{n+1,n+2}\circ\eta_{n,n+1}.
$$

The following straight sequence arises
$$
X^{\underrightarrow{\eta_{0, 1}}} I(X)\rightarrow ... \rightarrow
I^n(X)^{\underrightarrow{\eta_{n, n+1}}} I^{n+1}(X)\rightarrow
...\ . \eqno(1)
$$

Fix a metric $\rho$ on a compactum $X$ and the metrication
$\rho_{I,X}$ of the functor $I$. The metric on $I^n(X)$ generated
by this metrication denote through $\rho^n_{I,X}$. Then the maps
$$
\eta_{n,m}:(I^n(X),\rho^n_{I,X})\rightarrow(I^m(X),\rho^m_{I,X})
$$
are isometrical embeddings. The limit of the sequence (1) in
category metrizable spaces and their isometrical embeddings
denotes by $(I^+(X),\rho^+_{I,X})$. We give more constructive
definition of the metric $\rho^+_{I,X}$. By
$\eta_n:I^n(X)\rightarrow I^+(X)$ denotes the limit of embeddings
$\eta_{n,m}:I^n(X)\rightarrow I^m(X)$ under $m\rightarrow\infty$
consider while $I^+(X)$ as limit of the sequence (1) in the
category of sets. Then
$$
I^+(X)=\{\eta_n(I^n(X)):n\in\omega\},
$$
and the metric $\rho^+_{I,X}$ defines with metrics $\rho^n_{I,X}$
on the addends $\eta_n(I^n(X))$. More detail for
$x,y\in\eta_n(I^n(X))$ we have
$$
\rho^+_{I,X}(x,y)=\rho^n_{I,X}(a,b), \eqno(2)
$$
where $\eta_n(a)=x$, $\eta_n(b)=y$. The definition of the metric
$\rho^+_{I,X}$ through equality (2) is correct, since under $n<m$
the maps $\eta_{n,m}$ are isometrical embeddings.

If $f:X\rightarrow Y$ is a continuous then we can define the map
$I^+(f):I^+(X)\rightarrow I^+(Y)$. It does as the following way.
For $x\in I^+(X)$ there exists $n\in\omega$ and $a\in I^n(X)$ such
that $x=\eta_n(a)$. Put $I^+(f)(x)=\eta_n(I^n(X))(a)$. Since
$\eta_{n,m}$ is natural transformation of the functor $I^n$ into
the functor $I^m$ then this definition is correct.

Consider the following set
$$
I_f^{k+1}(X)=\{\mu\in I^{k+1}(X):\text{supp}\ \mu\subset
I_f^k(X),\ |\text{supp}\ \mu|<\omega\}.
$$

Analogously to linear case [2] idempotent probability measures
$\mu\in I_f^k(X)$ we call as measures with everywhere finite
supports. With recursion on $k$ it checks that $I_f^k(X)$ is
everywhere dense in $I^k(X)$.

\textbf{Lemma 4.} \textsl{Let $f:X\rightarrow Y$ be continuous
map, $k>0$. Then for all idempotent probability measures $^k\mu_1,
^k\mu_2\in I_f^k(X)$ the following inequality takes place}
$$
\rho^k_{\omega, Y}(I^k(f)(^k\mu_1),
I^k(f)(^k\mu_2))\leq\rho^k_{\omega, X}(^k\mu_1, ^k\mu_2).
$$

\textsc{Proof.} Let $^k\mu_1,\ ^k\mu_2\in I_f^k(X)$ be arbitrary
idempotent probability measures. Then there are $s_1,\ s_2\in N$
such that $supp(^k\mu_i)=\{^{k-1}\mu_{i1},...,
^{k-1}\mu_{is_i}\}$, $i=1,\ 2$, where $^{k-1}\mu_{il}\in
I^{k-1}(X)$, $l=1,\ ...,\ s_i$. Therefore the decompositions hold
$$
^k\mu_i=\lambda_{i1}\odot\delta_{^{k-1}\mu_{i1}} \oplus
…\oplus\lambda_{is_i}\odot\delta_{^{k-1}\mu_{is_i}},\,\ i=1,\ 2.
$$

According to the definition of the metric $\rho_I$ [1] we have
$$
\rho^k_{\omega, Y}(I^k(f)(^k\mu_1), I^k(f)(^k\mu_2))
\leq\rho^k_{\omega, X}(^k\mu_1, ^k\mu_2).
$$

Lemma 4 is proved.

Note, the inequality in Lemma 4 cannot replace  with equality.

\textbf{Example 1.} Let $X=Y=[0\ ,10]$, $\rho(t_1,
t_2)=|t_2-t_1|$, $t_1,\ t_2\in [0,1]$. Define the map
$f:X\rightarrow Y$ by formula
$$
f(x)=\left\{
\begin{array}{ll}
1-4\cdot\left(x-\frac{1}{2}\right)^2,& \mbox{ if }\  0\leq x\leq 1,\\
x-1, & \mbox{ if }\  1<x\leq 10.
\end{array}\right.
$$
We have
$$
f(0)=f(1)=0,\,\,\
f\left(\frac{1}{4}\right)=f\left(\frac{3}{4}\right)=f\left(\frac{7}{4}\right)=\frac{3}{4}.
$$

Define idempotent probability measures $\mu_1$ and $\mu_2$ by the
rules
$$
\mu_1=0\odot\delta_0\oplus(-5)\odot\delta_{\frac{1}{4}};\
\mu_2=0\odot\delta_{\frac{3}{4}}\oplus(-4)\odot\delta_1.
$$
It is easy to see that $supp(\mu_1)=\left\{0,\
\frac{1}{4}\right\}$ è $supp(\mu_2)=\left\{\frac{3}{4},\
1\right\}$. Then for each $\lambda\leq-5$ the idempotent
probability measure
$$
\xi_{\mu_1,\ \mu_2}=0\odot\delta_{\left(0, \frac{3}{4}\right)}
\oplus (-4)\odot\delta_{(0, 1)}\oplus (-5)
\odot\delta_{\left(\frac{1}{4}, \frac{3}{4}\right)}\oplus
\lambda\odot\delta_{\left(\frac{1}{4}, 1\right)}
$$
is an element of the set $\Lambda(\mu_1,\mu_2)$ (see [1]) which
satisfies Lemma 1 from [1]. That is why we have
$$
\rho_{\omega, X}(\mu_1,\mu_2)=5\frac{1}{2}.
$$

For any $\varphi\in C(Y)$ we have
$$
I(f)(\mu_1)(\varphi)=\mu_1(\varphi\circ
f)=\left(0\odot\delta_0\oplus(-5)\odot\delta_{\frac{1}{4}}\right)(\varphi\circ
f)=
$$
$$
=0\odot\delta_0(\varphi\circ
f)\oplus(-5)\odot\delta_{\frac{1}{4}}(\varphi\circ f)=
0\odot\varphi(f(0))\oplus(-5)\odot\varphi\left(f\left(\frac{1}{4}\right)\right)=
$$
$$
=0\odot\varphi(0)\oplus(-5)\odot\varphi\left(\frac{3}{4}\right)=
0\odot\delta_0(\varphi)\oplus(-5)\odot\delta_{\frac{3}{4}}(\varphi)=
\left(0\odot\delta_0\oplus(-5)\odot\delta_{\frac{3}{4}}\right)(\varphi).
$$
Hence
$I(f)(\mu_1)=0\odot\delta_0\oplus(-5)\odot\delta_{\frac{3}{4}}$.

Analogously it may be shown that
$I(f)(\mu_2)=(-4)\odot\delta_0\oplus 0\odot\delta_{\frac{3}{4}}$.

Thus $supp(I(f)(\mu_1))=supp(I(f)(\mu_2))=\left\{0,
\frac{3}{4}\right\}$. Here for any $\lambda\leq -5$ the idempotent
probability measure
$$
\xi_{I(f)(\mu_1),\ I(f)(\mu_2)}=0\odot\delta_{\left(0,\
\frac{3}{4}\right)}\oplus (-4)\odot\delta_{\left(0,\
0\right)}\oplus
(-5)\odot\delta_{\left(\frac{3}{4},\frac{3}{4}\right)}
\oplus\lambda\odot\delta_{\left(\frac{3}{4},\ 0\right)}
$$
is such an element of $\Lambda(I(f)(\mu_1),\ I(f)(\mu_2))$ which
satisfies Lemma 1 from [1]. That's why
$$
\rho_{\omega, Y}(I(f)(\mu_1),\ I(f)(\mu_2))=5.
$$

Thus $\rho_{\omega, Y}(I(f)(\mu_1),\ I(f)(\mu_2))\neq\rho_{\omega,
X}(\mu_1,\mu_2)$.

\textbf{Proposition 1.} \textsl{Let $X$, $Y$ be metric compacta
and what's more $diamX=diamY$. If a map $f:X\rightarrow Y$ is
$(\varepsilon,\ \delta)$-uniform continuous then the map
$I^k(f):I^k(X)\rightarrow I^k(Y)$  is also $(\varepsilon,\
\delta)$-uniform continuous.}

\textsc{Proof.} According to definition of the metric $\rho_{I,X}$
it is enough to establish the statement for idempotent probability
measures with everywhere finite supports. Without loss of
generality we can assume $\delta<\varepsilon$. But then Lemma 4
ends the proof. Proposition 1 is proved.

Finally we can formulate our main result.

$\textbf{Theorem 1.}$ \textsl{The functor $I$ has the following
properties:}

Ð1) \textsl{Let $(X_1,\ \rho^1)$ and  $(X_2,\ \rho^2)$ be metric
compacta. If $diam(X_1,\ \rho^1)=diam(X_2,\ \rho^2)$ and $i:(X_1,\
\rho^1)\rightarrow (X_1,\ \rho^1)$ is isometrical embedding then
$I(i):(I(X_1),\ \rho^1_{I, X_1})\rightarrow (I(X_1),\ \rho^1_{I,
X_2})$ is also isometric embedding;}

Ð2) \textsl{For any metric compactum $(X,\ \rho)$ the embedding
$\delta_X:(X,\ \rho)\rightarrow (I(X),\ \rho_{I, X})$ is an
isometry; } Ð3) \textsl{For any metric compactum $X$, and for an
arbitrary metric $\rho$ on $X$ the equality $diam(X,\
\rho)=diam(I(X),\ \rho_{I, X})$ holds;}

Ð4) \textsl{Let $(X_1,\ \rho^1)$ and $(X_2,\ \rho^2)$ be metric
compacta with $diamX_1=diamX_2$. Then for any continuous mapping
$f:(X_1,\ \rho^1)\rightarrow (X_2,\ \rho^2)$ the map
$I^+(f):(I^+(X_1),\ \rho^1_{I^+, X_1})\rightarrow (I^+(X_2),\
\rho^2_{I^+, X_2})$ is uniform continuous.}

\end{document}